\documentclass[12pt,reqno, final]{amsart}

\usepackage{a4}                          
\usepackage{a4wide}	
\usepackage[english]{babel}              
\usepackage[latin1]{inputenc}            
\usepackage{amsmath}                     
\usepackage{amssymb}                     
\usepackage{amsthm}
\usepackage{fancybox}
\usepackage{graphicx}
\usepackage{hyperref}
\usepackage{mathrsfs}									  
\usepackage[all]{xy}										
\usepackage[T1]{fontenc}								
\usepackage{url}												
\usepackage{array}
\usepackage{verbatim,epsfig}
\usepackage{wasysym}                    

\theoremstyle{plain}

\newtheorem{thm}{Theorem}[section] 
\newtheorem{cor}[thm]{Corollary}
\newtheorem{lem}[thm]{Lemma}

\theoremstyle{definition}

\newtheorem{rem}[thm]{Remark}

\newtheorem{question}[thm]{Question}

\newcommand{\NN}{{\mathbb N}}

\newcommand{\CC}{{\mathbb C}}

\newcommand{\MM}{{\mathbb M}}
\newcommand{\GG}{{\mathbb G}}
\newcommand{\HH}{{\mathbb H}}
\newcommand{\FF}{{\mathbb F}}

\renewcommand{\L}{{\mathscr L}}

\newcommand{\varps}{{\varepsilon}}

\newcommand{\im}{{\operatorname{im}}}
\newcommand{\htens}{\bar{\otimes}}

\newcommand{\tens}{\otimes}

\newcommand{\To}{\longrightarrow}

\newcommand{\red}{{\operatorname{red}}}
\newcommand{\Tor}{\operatorname{Tor}}
\newcommand{\op}{{\operatorname{{op}}}}

\newcommand{\Hom}{\operatorname{Hom}}

\newcommand{\id}{\operatorname{id}}

\newcommand{\ev}{{\operatorname{ev}}}

\newcommand{\bet}{\beta^{(2)}}

\newcommand{\alg}{{\operatorname{alg}}}

\newcommand{\Pol}{{\operatorname{Pol}}}
\newcommand{\deff}{{\operatorname{def}}}

\newcommand{\vN}{{\operatorname{vN}}}
\newcommand{\vn}{\vN}
\makeatletter
\def\equalsfill{$\m@th\mathord=\mkern-7mu
\cleaders\hbox{$\!\mathord=\!$}\hfill
\mkern-7mu\mathord=$}
\makeatother   

\begin{document}
\title{An $L^2$-K{\"u}nneth formula for tracial algebras}
\author{David Kyed}
\address{David Kyed\\ 
Mathematisches Institut\\
Georg-August-Universit{\"a}t G{\"o}ttingen\\
Bunsenstraße 3-5\\
37073 G{\"o}ttingen\\
Germany}

\email{kyed@uni-math.gwdg.de}
\urladdr{www.uni-math.gwdg.de/kyed}

\begin{abstract}
We prove a K{\"u}nneth formula computing the Connes-Shlyakhtenko $L^2$-Betti numbers of the algebraic tensor product of two tracial $*$-algebras in terms of the $L^2$-Betti numbers of the two original algebras. As an application, we construct examples of non-finite, non-cocommutative, compact quantum groups with a non-vanishing first $L^2$-Betti number.
\end{abstract}

\maketitle

\section{Introduction}
The theory of $L^2$-invariants originates from the work of Atiyah \cite{atiyah} in the 70's and was further developed by Cheeger and Gromov \cite{CG} in the 80's and  by L{\"u}ck \cite{Luck97,Luck98,Luck98-2} in the late 90's. L{\"u}ck showed that the whole theory  could be reformulated (and substantially extended) in a purely algebraic setting by extending the notion of Murray-von Neumann dimension from the class of finitely generated, projective modules over a finite von Neumann algebra $M$ to the class of \emph{all} $M$-modules. One of the pleasant features of L{\"u}ck's approach is that the algebraic setup allows the usage of all the powerful tools of homological algebra; for instance the $L^2$-homology  of a discrete group $\Gamma$ can be written as
\[
H_n^{(2)}(\Gamma)=\Tor_n^{\CC\Gamma}(\L(\Gamma),\CC),
\]
where $\L(\Gamma)$ denotes the group von Neumann algebra. The $L^2$-Betti numbers of $\Gamma$ are then obtained by applying the extended dimension function to the $L^2$-homology; in symbols $\bet_n(\Gamma)=\dim_{\L(\Gamma)}H_n^{(2)}(\Gamma)$. These $L^2$-Betti numbers permit a K{\"u}nneth formula \cite[2.7]{CG}; i.e.~for two discrete groups $\Gamma$ and $\Lambda$ we have that 
\[
\bet_n(\Gamma\times \Lambda)=\sum_{k+l=n}\bet_k(\Gamma)\bet_l(\Lambda).
\] 
In the beginning of the present century, Connes and Shlyakhtenko \cite{CS} took the development of $L^2$-invariants a step further by defining $L^2$-homology and $L^2$-Betti numbers for any weakly dense $*$-subalgebra $A$ in a tracial von Neumann algebra $(M,\tau)$; these are denoted $H_n^{(2)}(A,\tau)$ and $\bet_n(A,\tau)$ respectively. The Connes-Shlyakhtenko $L^2$-Betti numbers generalize the classical $L^2$-Betti numbers for groups by means of the formula
\[
\bet_n(\CC\Gamma,\tau)=\bet_n(\Gamma),
\]
where $\tau$ is the natural trace on the group von Neumann algebra $\L(\Gamma)$. The aim of the present note is to prove a K{\"u}nneth formula for the Connes-Shlyakhtenko $L^2$-Betti numbers; i.e.~that for weakly dense $*$-subalgebras $A$ and $B$ of tracial von Neumann algebras $(M,\tau)$ and $(N,\rho)$ we have
\[
\bet_n(A\odot B,\tau\tens \rho)=\sum_{k+l=n}\bet_k(A,\tau)\bet_l(B,\rho).
\]

\vspace{0.3cm}
\paragraph{\emph{Notation.}}
Above, and in what follows, $\odot$ is used to denote algebraic tensor products which, unless specified otherwise, are assumed to be over the complex numbers. The symbol $\tens$ will be reserved to denote the minimal tensor product of $C^*$-algebras, while $\htens$ will be used to denote the tensor product in the category of von Neumann algebras as well as the tensor product in the category of Hilbert spaces. Moreover, for any algebra $A$ we denote by $A^\op$ the opposite algebra and by $A^\ev$ the enveloping algebra $A\odot A^\op$. For any von Neumann algebra $M$ we let $M^{\overline{\ev}}$ denote the completed tensor product $M\htens M^\op$.
\vspace{0.3cm}
\paragraph{\emph{Structure.}}
The rest of the paper is organized in the following way: Section \ref{dim-section} is devoted to prove some minor results concerning the extended dimension function used in the proof of the K{\"u}nneth formula (Theorem \ref{kunneth}) which is presented in section \ref{kunneth-section}. In the fourth and final section we show how the K{\"u}nneth formula can be used to manufacture non-trivial compact quantum groups with a non-vanishing first $L^2$-Betti number.

\section{A bit of dimension theory}\label{dim-section}
In the present section we prove a few minor results related to L{\"u}ck's generalized Murray-von Neumann dimension $\dim_M(-)$; this is a dimension function defined on the category of all (algebraic) modules over a finite von Neumann algebra $M$ taking values in the interval $[0,\infty]$. We shall not dwell on the definition and (surprisingly nice!) properties of this dimension function, but refer the reader to Chapter 6 in \cite{Luck02}. All three results in this section are derived, without much effort, from the work of L{\"u}ck, but since they are not made explicit in the literature we present them here for the convenience of the reader.\\

Throughout this section, $M$ will denote a finite von Neumann algebra with a specified normal, faithful, tracial state $\tau$ and all calculations of Murray-von Neumann dimensions of $M$-modules are implicitly assumed to be with respect to the trace $\tau$. To fix notation, we recall that if $X$ is a submodule of an $M$-module $Y$ the algebraic closure of $X$ (relative to $Y$) is defined as
\[
\overline{X}^\alg=\bigcap_{\stackrel{\varphi\in \Hom(Y,M)}{X\subseteq \ker(\varphi)}} \ker(\varphi).
\]
Moreover, the projective part $P(X)$ of a module $X$ is defined as $X/\overline{\{0\}}^\alg$ and in case $X$ is finitely generated, $P(X)$ is in fact a finitely generated projective module. In the sequel, we will denote by $\pi_X$ the natural surjection of $X$ onto $P(X)$. As is easily seen, $P(-)$ becomes an endo-functor on the category of $M$-modules and it possesses the following property.
\begin{lem}\label{lem1}
For a homomorphism $f\colon X\to Y$ between finitely generated $M$-modules we have $\dim_M\im(f)=\dim_M\im(P(f))$.
\end{lem}
\begin{proof}
Since $\im(P(f))=\im(\pi_X\circ f)=\pi_X(\im(f))$ we have a short exact sequence
\begin{align}\label{ses}
0\To \ker(\pi_X|_{\im(f)})\overset{\subseteq}{\To} \im(f)\overset{\pi_X}{\To} \im(P(f))\To 0.
\end{align}
The kernel $\ker(\pi_X|_{\im(f)})$ is contained in the zero-dimensional $M$-module $\ker(\pi_X)=\overline{\{0\}}^\alg$, and therefore itself zero dimensional, and applying  additivity of $\dim_M(-)$ to the equation (\ref{ses}) we get
\begin{align*}
\dim_M\im(P(f))=\dim_M\im(f)-\dim_M\ker(\pi_X|_{\im(f)})=\dim_M\im(f).
\end{align*}
\end{proof}
\begin{lem}\label{lem2}
Let $f\colon P\to Q$ be a homomorphism of finitely generated projective $M$-modules and consider the continuous extension  $f^{(2)}\colon L^2(P)\To L^2(Q)$ of $f$ between the Hilbert $M$-module completions of $P$ and $Q$. Then $\dim_M\im(f)=\dim_M\overline{\im(f^{(2)})}^{\|\cdot\|_2}$. 
\end{lem}
Details about the notion of $L^2$-completions of projective modules can be found in \cite{Luck97}.

\begin{proof}
By \cite{Luck02} Theorem 6.24, the functor $L^2(-)$ going from the category of finitely generated projective $M$-modules to the category of finitely generated Hilbert $M$-modules is weakly exact and dimension preserving, with weakly exact and dimension preserving inverse, so
\[
\dim_{M}\overline{\im(f^{(2)})}^{\|\cdot\|_2}=\dim_M\overline{\im(f)}^\alg=\dim_M\im(f),
\]
where the last identity follows from \cite{Luck02} Theorem 6.7.
\end{proof}
The above two lemmas are included in order to prove the following result which will be essential in the proof of the K{\"u}nneth formula. The claim in Lemma \ref{lem3} was nested in the proof of \cite{Luck02} Theorem 6.54, but in order to clarify the proof of Theorem \ref{kunneth}  we have extracted the result as a separate lemma and included a proof.
\begin{lem}\label{lem3}
Let $F=(F_*,f_*)$ and $G=(G_*,g_*)$ be chain complexes consisting of finitely generated projective $M$-modules and consider a morphism of complexes $\varphi\colon F\to G$. Denote by $L^2(F)=(L^2(F_*),f_*^{(2)})$ and $L^2(G)=(L^2(G_*),g_*^{(2)})$ the completions of $F$ and $G$ into Hilbert $M$-chain complexes and by $\varphi^{(2)}\colon L^2(F)\to L^2(G)$ the morphism induced by $\varphi$. From this data we obtain three induced morphisms on the level on homology:
\begin{align*}
H_n(\varphi)&\colon H_n(F)\To H_n(G)\stackrel{\deff}{\hbox{\equalsfill}}\frac{\ker(g_n)}{\im(g_{n+1})};\\
\bar{H}_n(\varphi)&\colon \bar{H}_n(F)\To \bar{H}_n(G)\stackrel{\deff}{\hbox{\equalsfill}}\frac{\ker(g_n)}{\overline{\im(g_{n+1})}^\alg};\\
H_n^{(2)}(\varphi^{(2)})&\colon H_n^{(2)}(L^2(F))\To H_n^{(2)}(L^2(G))\stackrel{\deff}{\hbox{\equalsfill}}\frac{\ker(g_n^{(2)})}{\overline{\im(g_{n+1}^{(2)})}^{\|\cdot\|_2}}.
\end{align*}
The claim now is that $\dim_M\im(H_n(\varphi))=\dim_M\overline{\im(H_n^{(2)}(\varphi^{(2)}))}^{\|\cdot\|_2}$.
\end{lem} 

\begin{proof}
We first note that the homology modules $H_n(F)$ and $H_n(G)$ appearing in Lemma \ref{lem3} are finitely generated so that Lemma \ref{lem1} and Lemma \ref{lem2} applies; this is due to the fact that $M$ is a semihereditary ring and therefore \cite[0.2]{Luck97} its category of finitely presented modules is abelian. From \cite{Luck02} Lemma 6.52 we get an isomorphism $\bar{H}_n(F)\simeq P(H_n(F))$ under which $\bar{H}_n(\varphi)$ corresponds to $PH_n(\varphi)$ and an isomorphism  $L^2(PH_n(F))\simeq H_n^{(2)}(L^2(F))$ under which $(PH_n(\varphi))^{(2)}$ corresponds to $H_n^{(2)}(\varphi^{(2)})$. Hence
\begin{align*}
\dim_M\im(H_n(\varphi))&= \dim_M\im(PH_n(\varphi))\tag{Lemma \ref{lem1}}\\
&= \dim_M\overline{\im((PH_n(\varphi))^{(2)})}^{\|\cdot\|_2}\tag{Lemma \ref{lem2}}\\
&=\dim_M\overline{\im(H_n^{(2)}(\varphi^{(2)}))}^{\|\cdot\|_2}.
\end{align*}

\end{proof}

\section{The K{\"u}nneth formula}\label{kunneth-section}

Let $(M,\tau)$ and $(N,\rho)$ be tracial von Neumann algebras and let $A\subseteq M$ and $B\subseteq N$ be weakly dense $*$-subalgebras. Then the algebraic tensor product $A\odot B$ is a weakly dense $*$-subalgebra in the tracial von Neumann algebra $(M\htens N,\tau\tens \rho)$, and our aim now is to prove the following K{\"u}nneth formula for the Connes-Shlyakhtenko $L^2$-Betti numbers introduced in \cite{CS}.

\begin{thm}\label{kunneth}
For every $n\geq 0$ we have
\[
\bet_n(A\odot B,\tau\tens\rho)=\sum_{k+l=n}\bet_k(A,\tau)\bet_l(B,\rho),
\]
where $\bet_*(-,-)$ are the Connes-Shlyakhtenko $L^2$-Betti numbers of the tracial $*$-algebra in question. 
\end{thm}
Note that the $L^2$-Betti numbers might be infinite and the the above formula is therefore to be understood with respect to the standard rules for addition and multiplication in $[0,\infty]$.

\begin{proof}
Let $F=(F_*,f_*)$ and $G=(G_*,g_*)$ denote the bar-resolutions \cite[1.1.12]{loday} of $A$ and $B$ respectively and consider their tensor product $F\odot G=E$ which in degree $n$ has the module
\[
E_n =\bigoplus_{k+l=n} F_k\odot G_l,
\]
and whose $n$-th differential $e_n\colon E_{n}\to E_{n-1}$ is given by the formula
\[
e_n(x\tens y)=f_k(x)\tens y +(-1)^k x\tens g_l(y),
\]
for a homogeneous  element $x\tens y\in F_k\odot G_l$. Since both $F$ and $G$ are acyclic the same is true for $E$ (see e.g.~\cite[2.7.3]{weibel}) and $E$ therefore constitutes a resolution of $A\odot B$ in the category of $A\odot B$-bimodules --- a category we will freely identify  with the category of left modules over $(A\odot B)^{\ev}=(A\odot B)\odot (A\odot B)^\op$. As proven in \cite{CS} Lemma 2.2, the bar resolution $F$ can be written as an inductive limit of a family of subcomplexes $(F_{*,i},f_{*,i})_{i\in I}$ where each $F_i=(F_{*,i},f_{*,i})$ is a complex of finite length consisting of 
finitely generated free $A^\ev$-modules. If we denote the $n$-th homology of the complex
\[
(M\htens M^\op\odot_{A\odot A^\op} F_{i,*},1\tens f_{i,*})
\]
by $H_n(F_i^{\vn})$ and by $H_n(\varphi_{i_2i_1}^\vn)\colon H_n(F_{i_1}^{\vN})\to H_n(F_{i_2}^{\vN})$ the map induced by the inclusion $\varphi_{i_2i_1}:F_{i_1}\hookrightarrow F_{i_2}$ whenever $i_2\geq i_1$, then the $L^2$-homology $H_n^{(2)}(A,\tau)$  becomes the inductive limit
$\varinjlim(H_n(F_i^{\vN}),H_n(\varphi_{i_2i_1}^\vn))$. Since each $F_{n,i}$ is finitely generated, it follows from \cite{Luck02} Theorem 6.13 that
\begin{align}\label{beta-A}
\bet_n(A,\tau)&=\sup_{i_1}\inf_{\overset{\vspace{-0cm}}{i_2}\geq i_1}\dim_{M\htens M^\op}\im(H_n(\varphi_{i_2i_1}^\vn))\notag\\
&= \sup_{i_1}\inf_{\overset{\vspace{-0cm}}{i_2}\geq i_1}\dim_{M\htens M^\op}\overline{\im(H_n^{(2)}(\varphi_{i_2i_1}^{(2)}))}^{\|\cdot\|_2},
\end{align}
where the last equality follows from Lemma \ref{lem3}. In a completely similar manner, we obtain the $L^2$-homology $H_n^{(2)}(B,\rho)$ as the inductive limit $\varinjlim(H_n(G_j^{\vn}),H_n(\psi_{j_2j_1}^\vn))$ arising from a family of finite length subcomplexes $(G_{*,j},g_{*,j})_{j\in J}$, each of which consists of finitely generated free $B^\ev$-modules, and hence
\begin{align}\label{beta-B}
\bet_n(B,\rho)&=\sup_{j_1}\inf_{\overset{\vspace{-0cm}}{j_2}\geq j_1}\dim_{N\htens N^\op}\im(H_n(\psi_{j_2j_1}^\vn))\notag\\
&=\sup_{j_1}\inf_{\overset{\vspace{-0cm}}{j_2}\geq j_1}\dim_{N\htens N^\op}\overline{\im(H_n^{(2)}(\psi_{j_2j_1}^{(2)}))}^{\|\cdot\|_2}.
\end{align}
The families $(F_i)_{i\in I}$ and $(G_j)_{j\in J}$ defines a directed family of subcomplexes $(F_i\odot G_j )_{(i,j)\in I\times J}$ of $E$ \footnote{$I\times J$ is ordered by setting $(i_1,j_1)\geq (i_2,j_2)$ iff $i_1\geq i_2$ and $j_1\geq j_2$.} which has $E$ as its inductive limit. We now put $E_{ij}=F_i\odot G_j$ and denote by $H_n(E_{ij}^{\vn})$ the $n$-th homology of the induced complex
\[
\Big{(}(M\htens N)\htens(M\htens N)^\op\odot_{(A\odot B)^\ev} E_{ij,*}, 1\tens e_{ij,*}\Big{)}.
\]
For $(i_2,j_2)\geq (i_1,j_1)$, the inclusion $\varphi_{i_2i_1}\tens \psi_{j_2j_1}$ induces a map
\[
H_n((\varphi_{i_2i_1}\tens \psi_{j_2j_1})^\vn)\colon H_n(E_{i_1j_1}^{\vn})\To H_n(E_{i_2j_2}^{\vn}),
\]
and just as above we get that
\[
H_n^{(2)}(A\odot B,\tau\tens \rho)= \varinjlim\Big{(}H_n(E_{i,j}^{\vn}),H_n((\varphi_{i_2i_1}\tens \psi_{j_2j_1})^\vn)\Big{)}.
\]
Denoting the completed tensor product $(M\htens N)\htens (M\htens N)^\op$ by $(M\htens N)^{\overline{\ev}}$ we get
\begin{align}\label{beta-AB}
\bet_n(A\odot B,\tau\tens \rho)&=\sup_{(i_1,j_1)}\inf_{\overset{\vspace{0cm}}{(i_2,j_2)\geq (i_1,j_1)}}\dim_{(M\htens N)^{\overline{\ev}}}\im(H_n((\varphi_{i_2i_1}\tens \psi_{j_2j_1})^\vn))\notag\\
&= \sup_{(i_1,j_1)}\inf_{\overset{}{(i_2,j_2)\geq (i_1,j_1)}}\dim_{(M\htens N)^{\overline{\ev}}}\overline{\im(H_n^{(2)}((\varphi_{i_2i_1}\tens \psi_{j_2j_1})^{(2)}))}^{\|\cdot\|_2}.
\end{align}
For each $(i,j)\in I\times J$ we have two finitely generated Hilbert chain complexes; namely the Hilbert $(M\htens N)^{\overline{\ev}}$-chain complex $L^2(E_{ij})$  and the tensor product of Hilbert chain complexes $L^2(F_{i}) \htens L^2(F_j)$ (see e.g~\cite{Luck02}) which becomes a Hilbert chain complex for the  Neumann algebra $M^{\overline{\ev}}\htens N^{\overline{\ev}}$. The two von Neumann algebras in question are $*$-isomorphic in a trace-preserving way via the map
\[
(M\htens M^\op)\htens (N\htens N^\op)\ni a\tens c^\op\tens b\tens c^\op \overset{\alpha}{\longmapsto} a\tens b\tens c^\op\tens d^\op\in (M\htens N)\htens (M\htens N)^\op,
\]
and through the isomorphism $\alpha$ we  may therefore consider  $L^2(E_{ij})$ as a Hilbert $M^{\overline{\ev}}\htens N^{\overline{\ev}}$-chain complex; when doing so we write it as $_{\alpha}L^2(E_{ij})$. The Hilbert chain complex $_{\alpha}L^2(E_{ij})$   is nothing but the tensor product complex
$L^2(F_{i}) \htens L^2(G_j)$ and by Lemma 1.22 in \cite{Luck02} this identification gives rise to an $M^{\overline{\ev}}\htens N^{\overline{\ev}}$-isomorphism on the level of $L^2$-homology:
\[
_\alpha H_n^{(2)}(L^2(E_{ij}))=H_n^{(2)}(_\alpha L^2(E_{ij}))\overset{\sim}{\To}\bigoplus_{k+l=n} H^{(2)}_k(L^2(F_i))\htens H^{(2)}_l(L^2(G_j)).
\]
For each $(i_1,j_1)\leq (i_2,j_2)$ we therefore get a commutative diagram of Hilbert $M^{\overline{\ev}}\htens N^{\overline{\ev}}$-modules
\[
\xymatrix{
_\alpha H_n^{(2)}(L^2(E_{i_1j_1})) \ar^<<<<<<<<<{\sim}[rr]  \ar[dd]_{_\alpha H_n^{(2)}((\varphi_{i_2i_1}\tens \psi_{j_2j_1})^{(2)})} && \displaystyle\bigoplus_{k+l=n}^{\vspace{0.5cm}}H_k^{(2)}(L^2(F_{i_1}))\htens H_l^{(2)}(L^2(G_{j_1})) \ar[dd]^{\bigoplus_{k+l=n}^{\vspace{0.7cm}}H_k^{(2)}(\varphi_{i_2i_1}^{(2)})\tens H_l^{(2)}(\psi_{j_2j_1}^{(2)})     }        \\
            &&\\
_\alpha  H_n^{(2)}(L^2(E_{i_2j_2})) \ar^<<<<<<<<<{\sim}[rr] && \displaystyle\bigoplus_{k+l=n}H_k^{(2)}(L^2(F_{i_2}))\htens H_l^{(2)}(L^2(G_{j_2}))           
}
\]
For any finitely generated Hilbert $(M\htens N)^{\overline{\ev}}$-module $X$ we have
\[
\dim_{(M\htens N)^{\overline{\ev}}}X=\dim_{M^{\overline{\ev}}\htens N^{\overline{\ev}} }{_\alpha}X,
\]
simply because $\alpha$ is a trace-preserving $*$-isomorphism. We therefore get 
\begin{align}\label{sumformel}
&\dim_{(M\htens N)^{\overline{\ev}}}\overline{\im(H_n^{(2)}((\varphi_{i_2i_1}\tens \psi_{j_2j_1})^{(2)}))}^{\|\cdot\|_2}\notag\\
&=\dim_{M^{\overline{\ev}}\htens N^{\overline{\ev}}}\overline{\im(_\alpha H_n^{(2)}((\varphi_{i_2i_1}\tens \psi_{j_2j_1})^{(2)}))}^{\|\cdot\|_2}\notag\\
&= \dim_{M^{\overline{\ev}}\htens N^{\overline{\ev}}}\overline{\im(\oplus_{k+l=n}H_k^{(2)}(\varphi_{i_2i_1}^{(2)})\tens H_l^{(2)}(\psi_{j_2j_1}^{(2)}))}^{\|\cdot\|_2}\notag\\
&=\dim_{M^{\overline{\ev}}\htens N^{\overline{\ev}}}\bigoplus_{k+l=n}\overline{\im(H_k^{(2)}(\varphi_{i_2i_1}^{(2)}))}^{\|\cdot\|_2}\bar{\otimes}\hspace{0.1cm} \overline{\im(H_l^{(2)}(\psi_{j_2j_1}^{(2)}))}^{\|\cdot\|_2}\notag\\
&=\sum_{k+l=n}\dim_{M^{\overline{\ev}}}\Big{(}\overline{\im(H_k^{(2)}(\varphi_{i_2i_1}^{(2)}))}^{\|\cdot\|_2}\Big{)}
\dim_{N^{\overline{\ev}}}\Big{(}\overline{\im(H_l^{(2)}(\psi_{j_2j_1}^{(2)}))}^{\|\cdot\|_2}\Big{)},
\end{align}
where the last equality follows from \cite{Luck02} Theorem 1.12. Combining all the formulas obtained so far, the desired identity follows:
\begin{align*}
&\bet_n(A\odot B,\tau\tens \rho)=\\
&\overset{(\ref{beta-AB})}{=} \sup_{(i_1,j_1)}\inf_{\overset{}{(i_2,j_2)\geq (i_1,j_1)}}\dim_{(M\htens N)^{\overline{\ev}}}\overline{\im(H_n^{(2)}((\varphi_{i_2i_1}\tens \psi_{j_2j_1})^{(2)}))}^{\|\cdot\|_2}\\
&\overset{(\ref{sumformel})}{=}\sup_{(i_1,j_1)}\inf_{\overset{}{(i_2,j_2)\geq (i_1,j_1)}} \sum_{k+l=n}\dim_{M^{\overline{\ev}}}\Big{(}\overline{\im(H_k^{(2)}(\varphi_{i_2i_1}^{(2)}))}^{\|\cdot\|_2}\Big{)}
\dim_{N^{\overline{\ev}}}\Big{(}\overline{\im(H_l^{(2)}(\psi_{j_2j_1}^{(2)}))}^{\|\cdot\|_2}\Big{)}\\
&=\sum_{k+l=n}\Big{(}\sup_{i_1}\inf_{\overset{}{i_2\geq i_1}}\dim_{M^{\overline{\ev}}}\overline{\im(H_k^{(2)}(\varphi_{i_2i_1}^{(2)}))}^{\|\cdot\|_2}\Big{)}
\Big{(}\sup_{j_1}\inf_{\overset{}{j_2\geq j_1}}\dim_{N^{\overline{\ev}}}\overline{\im(H_l^{(2)}(\psi_{j_2j_1}^{(2)}))}^{\|\cdot\|_2}\Big{)}\\
&\overset{(\ref{beta-A},\ref{beta-B})}{=}\sum_{k+l=n}\bet_k(A,\tau)\bet_l(B,\rho)
\end{align*}
\end{proof}
The K{\"u}nneth formula gives an alternative proof of the following well known fact which is usually derived from Theorem 2.4 in \cite{CS}.
\begin{cor}
The $n$-th $L^2$-Betti number of the hyperfinite factor $R$ is either zero or infinite. 
\end{cor}
\begin{proof}
Denote by $\tau$ the unique trace on $R$ and by $\rho$ the normalized trace on $\MM_2(\CC)$. Since $R$ is hyperfinite it absorbs $\MM_2(\CC)$ tensorially and since both $R$ and $\MM_2(\CC)$ are factors any isomorphism $R\simeq R\odot \MM_2(\CC)$ is bound to preserve the traces. Moreover, it follows from \cite{CS} Proposition 2.9 that
\[
\bet_k(\MM_2(\CC),\rho)=
\begin{cases}
\frac14 &\text{ if }k=0, \\
0 &\text{ if }k\geq 1. 
\end{cases}
\]
Hence
\[
\bet_n(R,\tau)=\bet_n(R\odot \MM_2(\CC),\tau\tens\rho)=\sum_{k+l=n}\bet_k(R,\tau)\bet_l(\MM_2(\CC),\rho)=\bet_n(R,\tau)\tfrac{1}{4},
\]
which forces $\bet_n(R,\tau)$ to be either zero or infinite.
\end{proof}
It is, to the best of the author's knowledge, still not known what the $L^2$-Betti numbers of the hyperfinite factor are, except in degree zero where it follows from \cite{CS} Corollary 2.8 that $\bet_0(R,\tau)=0$. However, having in mind the well known analogy between hyperfiniteness and amenability, it is of course natural to expect that also the higher $L^2$-Betti numbers of $R$ vanish.

\section{An application towards quantum groups}
We take as our starting point Woronowicz's definition \cite{wor-cp-qgrps} of a compact quantum group. Thus, a compact quantum group $\GG$ consists of a (not necessarily commutative) unital, separable $C^*$-algebra $C(\GG)$ together with a coassociative, unital $*$-homomorphism $\Delta_\GG\colon C(\GG)\to C(\GG)\tens C(\GG)$, called the comultiplication, which furthermore has to satisfy a certain non-degeneracy condition. We remind the reader that such a $C^*$-algebraic quantum group automatically gives rise to a purely algebraic quantum group (a Hopf $*$-algebra \cite{klimyk}), whose underlying algebra will be denoted $\Pol(\GG)$, as well as a von Neumann algebraic quantum group \cite{kustermans-vaes} whose underlying algebra will be denoted $L^\infty(\GG)$. We also recall that the $C^*$-algebra $C(\GG)$ comes with a distinguished state $h_\GG$, called the Haar state, which plays the role corresponding to the Haar measure on a genuine, compact group. Performing the GNS construction with respect to the Haar state yields a Hilbert space denoted $L^2(\GG)$. The canonical example of a compact quantum group, on which the general definition is modeled,  is obtained by considering a compact, second countable, Hausdorff topological group $G$ and its commutative $C^*$-algebra $C(G)$ of continuous, complex valued functions. In this case, the von Neumann algebra becomes $L^\infty(G,\mu)$, where $\mu$ denotes the Haar probability measure, and the associated Hopf $*$-algebra becomes the algebra generated by matrix coefficients arising from the irreducible representations of $G$. If $\GG$ and $\HH$ are two compact quantum groups they give rise to a third quantum group, denoted $\GG\times \HH$, whose underlying $C^*$-algebra is $C(\GG)\tens C(\HH)$ and whose comultiplication is given by $\Delta_{\GG\times \HH}=(\id\tens\sigma\tens \id)\Delta_\GG\tens \Delta_\HH$; here $\sigma\colon C(\GG)\tens C(\HH)\to C(\HH)\tens C(\GG)$ denotes the flip-isomorphism. The Haar state $h_{\GG\times \HH}$ is given by $h_\GG\tens h_\HH$ and $\Pol(\GG\times \HH)=\Pol(\GG)\odot \Pol(\HH)$. See \cite{wang-tensor} for more details.\\

In \cite{quantum-betti} and \cite{coamenable-betti} the notion of $L^2$-invariants was generalized from the class of discrete groups to the class of compact quantum groups of Kac type; i.e.~those whose Haar state is a trace. If $\GG$ is such a quantum group its $n$-th  $L^2$-homology $H^{(2)}_n(\GG)$ is defined as $\Tor_n^{\Pol(\GG)}(L^\infty(\GG),\CC)$ and the $n$-th $L^2$-Betti number $\bet_n(\GG)$ is defined by applying L{\"u}ck's generalized Murray-von Neumann dimension $\dim_{L^\infty(\GG)}(-)$ to the $L^\infty(\GG)$-module $H_n^{(2)}(\GG)$. As a consequence of Theorem \ref{kunneth}, we also obtain a K{\"u}nneth formula for these quantum group $L^2$-Betti numbers.

\begin{cor}\label{quantum-kunneth}
Let $\GG$ and $\HH$ be compact quantum groups of Kac type. Then 
\[
\bet_n(\GG\times \HH)=\sum_{k+l=n}\bet_k(\GG)\bet_l(\HH),
\]
for every $n\geq 0$.
\end{cor}
\begin{proof}
It was shown in \cite{quantum-betti} Theorem 4.1 that $\bet_n(\GG)$ coincides with the $n$-th Connes-Shlyakh\-tenko $L^2$-Betti number $\bet_n(\Pol(\GG),h_\GG)$, and since $\Pol(\GG\times \HH)=\Pol(\GG)\odot \Pol(\HH)$ and $h_{\GG\times \HH}=h_{\GG}\tens h_{\HH}$ \cite{wang-tensor}  the result follows from Theorem \ref{kunneth}.
\end{proof}

Next we explain how the K{\"u}nneth formula provides us with non-trivial examples of quantum groups with non-vanishing $L^2$-Betti numbers.\\

	All cocommutative quantum groups are (in their reduced form) isomorphic to a quantum group of the form $(C^*_\red(\Gamma),\Delta_\red)$ where $\Gamma$ is a discrete group and $\Delta_\red(\gamma)=\gamma\tens\gamma$. It follows from the definitions \cite[1.3]{quantum-betti} that the $L^2$-Betti numbers of such a quantum group coincides with the classical $L^2$-Betti numbers of the underlying group $\Gamma$. Considering, for instance, the case of the free group on two generators $\FF_2$ with $\bet_1(\FF_2)=1$ we therefore get, in a somewhat trivial way, a compact quantum group with a non-vanishing first $L^2$-Betti number. Another trivial source of non-vanishing results is the class of finite quantum groups (i.e.~those whose $C^*$-algebra is finite dimensional); for such a quantum group $\GG$ the zeroth $L^2$-Betti number equals $\dim_\CC(C(\GG))^{-1}$ and all the higher $L^2$-Betti numbers vanish \cite{quantum-betti}.  So far, these are the only known examples of quantum groups with non-vanishing $L^2$-Betti numbers and the following question is therefore natural. 

\begin{question}
What is an example of a non-finite, non-cocommutative compact quantum group with a positive $L^2$-Betti number? 
\end{question}
The K{\"u}nneth formula provides an answer to this question. For this, let $\GG$ be a finite, non-cocommutative quantum group and denote by $N$ the dimension of $C(\GG)$ and let $\HH$ be the compact quantum group arising from $\FF_2$. We have
\[
\bet_p(\GG)=\left\{%
\begin{array}{ll}
   \frac{1}{N}, & \hbox{when $p=0$;} \\
    0, & \hbox{otherwise;} \\
\end{array}%
\right. \qquad \text{ and } \qquad  \bet_p(\HH)=\left\{%
\begin{array}{ll}
   1, & \hbox{when $p=1$;} \\
    0, & \hbox{otherwise,} \\
\end{array}%
\right.
\]
and the K{\"u}nneth formula therefore yields
\[
\bet_1(\GG\times \HH)=\bet_0(\GG)\bet_1(\HH) +\bet_1(\GG)\bet_0(\HH)=\tfrac{1}{N}.
\]
By construction, $C(\GG\times \HH)$ has infinite linear dimension and since $\GG$ is assumed non-cocommutative $\GG\times \HH$ becomes non-cocommutative.
\begin{rem}
For the  free group on  $k$ generators $\FF_k$ the only non-vanishing $L^2$-Betti number is the first which has value  $k-1$. Also, for each $n\in \NN$ it is easy to produce a finite quantum group of dimension $n$; one may simply take a group $G$ with $n$ elements and consider the associated commutative quantum group $C(G)$. By copying the example from  above we may therefore construct quantum groups  with any prescribed positive, rational number as its first $L^2$-Betti number. Note, however, that if  the group $G$ is chosen (or forced) to be abelian the example becomes cocommutative.
\end{rem}
As another consequence of the K{\"u}nneth formula we also get the following vanishing result.

\begin{cor}
Let $\GG$ and $\HH$ be compact quantum groups of Kac type and assume one of them to be infinite and coamenable. Then $\bet_n(\GG\times \HH)=0$ for all $n\geq 0$.
\end{cor}
We remind the  reader that a compact quantum group is called coamenable \cite{murphy-tuset} if the counit $\varps\colon \Pol(\GG)\to\CC$ extends to a character on the image under the GNS-representation of $C(\GG)$ on $L^2(\GG)$. 

\begin{proof}
Corollary 6.2 in \cite{coamenable-betti} together with \cite{my-thesis} Proposition 5.1.5 shows that all $L^2$-Betti numbers of an infinite, coamenable quantum group vanish and the claim now follows from Corollary \ref{quantum-kunneth}.

\end{proof}

\end{document}